\input amstex
\documentstyle{amsppt}

%
\newcount\eqnumber
\eqnumber=1

\def\neweq{{\rm{\the\eqnumber}}
\global\advance\eqnumber by 1}

\def\eqnam#1{\xdef#1{\the\eqnumber}}


\topmatter
\title {A pointwise estimate for the Fourier transform and maxima of a function}\endtitle
\leftheadtext{R. Berndt}
\rightheadtext{Fourier transform and Maxima}
\author {Ryan Berndt}\endauthor

\address {Yale University and Otterbein College
}\endaddress
\curraddr {Otterbein College, Westerville, Ohio 43081}
\endcurraddr 
\email {rberndt\@otterbein.edu}\endemail




\dedicatory \enddedicatory

\abstract We show a pointwise estimate for the Fourier transform on the line involving the number of times the function changes monotonicity. The contrapositive of the theorem may be used to find a lower bound to the number of local maxima of a function. We also show two applications of the theorem. The first is the two weight problem for the Fourier transform, and the second is estimating the number of roots of the derivative of a function.
\endabstract

\endtopmatter

\document

It is a classical result of Dirichlet  that if $f$ is a function of bounded variation on the circle, then the Fourier coefficients, $\widehat f(n)$, are $O(1/n)$ (and therefore the Fourier series of $f$ converges) \cite{9, p. 128}, \cite{10, p. 57}. We present here an inequality that implies a similar result, but for the Fourier transform on the line. Each time a real function changes from increasing to decreasing, we say that the function {\it crests}. We show an estimate for the Fourier transform of a function in terms of the number times the function crests.

This paper consists of two theorems and two applications. The first application estimates the number of roots of the derivative of a function, and the second application is a weighted Fourier norm inequality.

We first make a quick note on terminology and notation. We use the terms increasing and decreasing in the wider sense; $f(x)\equiv 1$ is both increasing and decreasing everywhere. We define the Fourier transform by the formula $\widehat f(z)=\int f(x)e^{-ixz}\,dx$. Whenever we take the Fourier transform of a function, we assume that $f\in L^1$ so that $\widehat f(z)$ is defined for all $z\in \bold{R}$. We use the letter $C$ to denote a constant whose value may change at each appearance. Finally, we say that two sets have almost disjoint support if the intersection of their supports has Lebesgue measure zero.

We provide a precise definition of crests below, but the reader may want to think of them as local maxima for the time being.

\proclaim{Theorem 1} If $f\in L^1$ is nonnegative and $\#crests(f)  = N$ then

$$|\widehat f(z)|\leq N\pi\sqrt{10}\int_0^{1/z} f^*(x)\,dx.$$

\noindent for all $z>0$.
\endproclaim

Here, $f^*$ is the decreasing rearrangement of $f$. As usual, it is defined by 
$f^*(x)=\inf \{\alpha >0 :|\{t:|f(t)|>\alpha\}|\leq x\}$, where $|\cdot|$ of a set represents the Lebesgue measure of that set. We note that if $f$ is also bounded then the theorem implies that $\widehat f(z)$ is $O(1/z)$, as in the case of Fourier series \cite{9, p. 128}.

In an example below we demonstrate that the appearance of the $N$ in the theorem can not be removed, and in fact, appears as the correct order of magnitude. Therefore, we are able to turn our viewpoint and use the contrapositive to predict the number of times that that the function will crest. Precisely, the contrapositive is the following.

\proclaim{Theorem 2} If $f\in L^1$ is nonnegative and 

$$Q(z) =  \frac{|\widehat f(z)|}{\pi\sqrt{10}\int_0^{1/z} f^*(x)\,dx}>N$$

\noindent for some $z>0$, then $\#crests(f) > N$.
\endproclaim

We prove this theorem below. We note that the function $Q$ is continuous since $\widehat f$ is continuous and the integral is absolutely continuous. So, if $Q(z)>N$ for some $z$ then it is greater than $N$ in a neighborhood of $z$. Application 2, below, shows how we may use Theorem 2 to estimate the number of roots of the derivative of a function $f$.

\definition{Definition}  A nonnegative function $f$ is said to {\it crest once} if there exists a point $b$ such that $f(x)$ is increasing for $x < b$ and decreasing for $x > b$. In this case we write $\#crests(f)=1$.
\enddefinition

\definition{Definition} We say that a nonnegative function $f$ {\it crests} $N \geq 1$ {\it times} if it can be written as no fewer than the sum of $N$ nonnegative functions with almost disjoint support, each of which crests once. That is,

$$\#crests(f)=\min\{N\in \bold{N}:f=\sum_{i=1}^N f_i, f_i\geq 0 \hbox{\ al. disj. supp.}, \#crests(f_i)=1\}
$$

\noindent If the set above is empty then we say that $\#crests(f)=\infty$.

\enddefinition

The sum of two disjoint characteristic functions like

$$f(x)=\chi_{[0,1]}(x)+\chi_{[2,3]}(x)$$

\noindent crests two times. If $f$ is zero on the negative axis and decreasing as $x$ grows then $f$ crests once. For example,

$$f(x) = \cases 0 & \hbox{for\ } x\leq 0\cr
1/x& \hbox{for\ }x> 0.
\endcases$$

\noindent has one crest.  A constant function has one crest, and if $f$ is a strictly increasing function, $\#crests(f)=\infty$.

Sometimes the number of crests equals the number of local maxima of a function. Any condition on a function that forces it to be locally strictly increasing and decreasing near a maximum will imply that the number of crests equals the number of local maxima. For example, if $f$ is a smooth function such that $f'(x)=0$ implies $f''(x)\ne 0$, then $\#crests(f)$ equals the number of local maxima of $f$.\medskip

\demo{Proof of Theorem 2} This is really just the contrapositive of Theorem 1. Suppose $f\in L^1$ is nonnegative and $Q(z)>N$ for some $z>0$, then by Theorem 1, $\#crests(f)\ne N$. Either $\#crests(f)>N$ or $\#crests(f)<N$. If $\#crests(f)<N$ then it must be possible to write $f$ as the sum of fewer than $N$ functions, each with one crest, with almost disjoint supports. But, then, by the theorem, $Q(z)<N$ for all $z>0$, a contradiction. Hence, $\#crests(f)>N$.

\enddemo

We prove Theorem 1 by first proving two lemmata. We start by considering the case where $f$ is a decreasing function and use this to bootstrap to the case of a finite number of crests. We note that by $L^1[0,\infty)$ we mean the space of all integrable functions that are zero on the negative axis.

\proclaim{Lemma} If $f\in L^1[0,\infty)$ is nonnegative and decreasing then

$$|\widehat f(z)| \leq \frac{\pi}{2}\sqrt{10}\int_0^{1/z} f(x)\,dx\tag \eqnam\basicfourier\neweq$$

\noindent for all $z>0$.

\endproclaim

\demo{Proof} Since $f$ is zero on the negative axis we may express the Fourier transform as the difference of the cosine and sine transforms:

$$\aligned \widehat f(z) &= \int_0^\infty f(x)\cos(xz)\,dx - i \int_0^\infty f(x)\sin(xz)\,dx\cr
&=Cf(z)-iSf(z).\cr
\endaligned$$

\noindent We prove (\basicfourier) by showing 

$$0< Sf(z) \leq \int_0^{\pi/2z} f(x)\,dx\hbox{\ \ and\ \ } |Cf(z)|\leq \int_0^{3\pi/2z} f(x)\,dx\tag\eqnam\basicsine\neweq$$

\noindent for all $z>0$. Since $f$ is decreasing (\basicsine) implies $Sf(z)\leq \pi/2 \int_0^{1/z} f(x)\,dx$ and $|Cf(z)| \leq 3\pi/2\int_0^{1/z} f(x)\,dx$. Therefore, $|\widehat f(z)| =\sqrt{Sf(z)^2+Cf(z)^2} \leq \break \sqrt{\pi^2/4+9\pi^2/4}\int_0^{1/z} f(x)\,dx =  \pi\sqrt{10}/2\int_0^{1/z}f(x)\,dx$.

To prove the inequality (\basicsine) for the sine transform, we fix $z>0$ and write it as an alternating series in the following way:

$$\aligned
Sf(z) &= \int_0^\infty f(x)\sin(xz)\,dz\cr 
&= \frac{1}{z}\int_0^\infty f(x/z)\sin x \,dx\cr
&=\frac{1}{z}\sum_{k=0}^\infty (-1)^k\int_{k\pi}^{(k+1)\pi} f(x/z)|\sin x |\,dx\cr
&=\frac{1}{z}\sum_{k=0}^\infty (-1)^k b_k,\cr
\endaligned$$

\noindent where $b_k=\int_{k\pi}^{(k+1)\pi} f(x/z)|\sin(x)|\,dx \geq 0$. Since $f$ is decreasing, $b_k$ is a decreasing sequence. Therefore, by a standard alternating series estimate,

$$0 \leq  b_0-b_1 \leq zSf(z) \leq b_0,$$

\noindent proving (\basicsine). 

The same technique is used to prove the estimate for the cosine transform. Let $\alpha_k=(k+1/2)\pi$. Fixing $z>0$, we write

$$\aligned
Cf(z) &=\int_0^\infty f(x)\cos(xz)\,dx\cr
&=\frac{1}{z}\int_0^\infty f(x/z)\cos x \,dx\cr
&=\frac{1}{z}\int_0^{\pi/2} f(x/z)\cos x \,dx +\frac{1}{z}\sum_{k=1}^\infty (-1)^k\int_{\alpha_k}^{\alpha_k+\pi} f(x/z)|\cos x |\,dx\cr
&=\frac{a_0}{z} +\frac{1}{z}\sum_{k=1}^\infty (-1)^k a_k,\cr
\endaligned$$

\noindent where $a_0=\int_0^{\pi/2} f(x/z)\cos(x)\,dx$ and $a_k=\int_{\alpha_k}^{\alpha_k+\pi} f(x/z)\cos(x)\,dx$ for $a_k\geq 1$. Since $f$ is decreasing, we know that $a_1 \geq a_2 \geq a_3 \geq \dots$. Because the intervals over which they are defined are of different lengths, we do not know how $a_0$ compares to $a_k$, for $k\geq 1$. But, 

$$a_0 - zCf(z) = a_1-a_2+a_3-\dots$$

\noindent is an alternating series on which we can apply the same standard estimate as before to get

$$0\leq a_0-zCf(z) \leq a_1.$$

\noindent Since $a_0 \geq 0$, we have that $-a_1\leq zCf(z)\leq a_0$, and thus $|Cf(z)| \leq (a_0+a_1)/z$, finishing the proof of (\basicsine), and therefore the lemma.
\enddemo

\proclaim{Lemma} If $f\in L^1[0,\infty)$ is nonnegative and crests once at $x=b$ then

$$|\widehat f(z)| \leq \frac{\pi}{2}\sqrt{10}\int_{b-1/z}^{b+1/z} f(x)\,dx\tag \eqnam\onepeak\neweq$$

\noindent for all $z>0$.
\endproclaim

\demo{Proof} We may write $f=g_1+g_2$ where $g_1$ is supported in $[a,b]$ and increasing over its support, and $g_2$ is supported in $[b,\infty)$ and decreasing over its support. If we let $h(x)=g_1(b-x)$ then $h$ is decreasing and we may apply (\basicfourier) to $h$ to get

$$\aligned
|\widehat h(z)| &\leq \frac{\pi}{2}\sqrt{10}\int_0^{1/z} h(x)\,dx\cr
&=\frac{\pi}{2}\sqrt{10}\int_{b-1/z}^b g_1(x)\,dx.\cr
\endaligned$$

\noindent Since $\widehat h(z) = e^{-ibz}\widehat g_1(-z)$, we have $|\widehat h(z)| = |\widehat g_1(-z)| = |g_1(z)|$. Hence,

$$|\widehat g_1(z)| \leq \frac{\pi}{2}\sqrt{10}\int_{b-1/z}^b g_1(x)\,dx.$$

\noindent Similarly, we let $h(x)=g_2(x+b)$. Then, $h$ is decreasing and we may apply (\basicfourier) and the fact that $|\widehat h(z)|=|\widehat g_2(z)|$ to get

$$|\widehat g_2(z)| \leq \frac{\pi}{2}\sqrt{10}\int_b^{b+1/z} g_2(x)\,dx.$$

\noindent We apply the triangle inequality to $\widehat f=\widehat g_1+\widehat g_2$ to finish the proof.
\enddemo

\demo{Proof of Theorem 1} We define 

$$f_1(x) = \cases 0, & \hbox{for\ } x< 0;\cr
f(x) & \hbox{for\ }x\geq 0;\cr
\endcases$$

$$f_2(x) = \cases 0, & \hbox{for\ } x< 0;\cr
f(-x) & \hbox{for\ }x\geq 0,\cr
\endcases$$

\noindent so that $f(x)=f_1(x)+f_2(-x)$ and $f_1,f_2 \in L^1[0,\infty)$. Let $N_j=\#crests(f_j)$. Since the supports of $f_j$ overlap only at $x=0$, $N=N_1+N_2$. Also, there exist functions $f_{j,i}$ such that $0\leq f_{j,i}(x)\leq f_j(x) \leq f(x)$, $\#crests(f_{i,j})=1$, and

$$f_j(x) = \sum_{i=1}^{N_j} f_{j,i}(x).$$

\noindent Applying the linearity of the Fourier transform and the fact that the modulus of the Fourier transform of a real function is even, we have $|\widehat f(z)| = |\widehat f_1(z)+\widehat f_2(-z)| \leq |\widehat f_1(z)|+|\widehat f_2(z)|$. Supposing that the functions $f_{1,i}$ and $f_{2,i}$ have $b_i$ and $c_i$, respectively, as their cresting points, we have

$$\aligned
|\widehat f(z)| 
&\leq \frac{\pi}{2}\sqrt{10}\left(\sum_{i=1}^{N_1} \int_{b_i-1/z}^{b_i+1/z} f_{1,i}(x)\,dx+ \sum_{i=1}^{N_2} \int_{c_i-1/z}^{c_i+1/z} f_{2,i}(x)\,dx\right)\cr
\endaligned$$

\noindent with the help of repeated applications of (\onepeak). For any Lebesgue measurable set $E$, $\int_E f \leq \int_0^{|E|} f^*$, see Bennett and Sharpley \cite{2, p. 44}. Thus,

$$\aligned |\widehat f(z)|
&\leq \frac{\pi}{2}\sqrt{10}\left(\sum_{i=1}^{N_1} \int_{b_i-1/z}^{b_i+1/z} f(x)\,dx+ \sum_{i=1}^{N_2} \int_{c_i-1/z}^{c_i+1/z} f(x)\,dx\right)\cr
&\leq \frac{\pi}{2}\sqrt{10}\left(\sum_{i=1}^{N_1} \int_{0}^{2/z} f^*(x)\,dx+ \sum_{i=1}^{N_2} \int_{0}^{2/z} f^*(x)\,dx\right)\cr
&\leq N\pi\sqrt{10}\int_0^{1/z} f^*(x)\,dx,\cr
\endaligned$$
\enddemo

\noindent finishing the proof of the Theorem 1.

\example{Example} In this example we show that the $N$ in Theorem 1 can not be removed and appears as the correct order of magnitude. Precisely, we show that given $N\geq 1$, there exists a function $f\in L^1[0,\infty)$ with $5N$ crests such that

$$Q(z) = \frac{|\widehat f(z)|}{\pi\sqrt{10}\int_0^{1/z} f^*(x)\,dx} > N$$

\noindent for some $z>0$. We take $f$ to be of the form

$$f(x)=\sum_{k=0}^\infty c_k\chi_{[k,k+1]}(x).$$

\noindent We let $c_0, c_2, c_4, \dots, c_{2(5N-1)}$ be $1$ and we let all other $c_k$ be zero.
The Fourier transform of our function is given by

$$\widehat f(z) = \frac{1}{z}\sum_{k=0}^\infty c_k[\sin(kz+z)-\sin(kz)]-ic_k[\cos(kz)-\cos(kz+z)],$$

\noindent and the decreasing rearrangement is given by

$$f^*(x) = \cases 1 & \hbox{for\ } x< 5N\cr
0 & \hbox{for\ }x \geq 5N.
\endcases$$

\noindent Now, $f$ is a function with $5N$ crests, but if we take $z=2l\pi$, $l\in \bold{N}$ then

$$\widehat f(z) = \frac{-2i}{z}(c_0+c_2+c_4+\dots) = \frac{-10Ni}{z},$$

\noindent and if $z$ is also greater than $1/5N$

$$\int_0^{1/z} f^*(x)\,dx =1/z.$$

\noindent Hence, for large enough $z=2l\pi$,

$$Q(z)=\frac{10N/z}{\pi\sqrt{10}/z} \approx 1.007N >N.$$
\endexample

\example{Application 1} In view of this example, we can use Theorem 2 to estimate the number of roots of the derivative of a function.

Suppose $f$ is a smooth function where $f'(x)=0$ implies $f''(x)\ne 0$, so that the derivative crosses the $x$-axis at each of its roots. In this case, the number of crests is equal to the number of local maxima of $f$. Now, if $f$ is integrable and has $N$ local maxima, then $f$ has at least $2N-1$ local extrema and $f'$ has at least $2N-1$ roots. Hence, we may formulate the following application of our theorem.

\proclaim{Corollary} Suppose $f\in L^1$ is nonnegative, smooth, and $f'(x)=0$ implies $f''(x)\ne 0$. If $Q(z)>N$ for some $z>0$ then $f'$ has at least $2N-1$ real roots.
\endproclaim
\endexample

\example{Application 2} In this application we show how we can apply the heart of Theorem 1, appearing in inequality (\basicfourier), to a norm estimate for the Fourier transform. The norm estimate we have in mind is the ``two weight problem for the Fourier transform." Part of this problem is finding functions $u$ and $v$ and a constant $C$ such that

$$\left(\int |\widehat f(z)|^q u(z)\,dz\right)^{1/q} \leq C \left(\int |f(x)|^pv(x)\,dx\right)^{1/p}\tag \eqnam\twoweight\neweq$$

\noindent for all $f$ where the right hand side is finite and the Fourier transform is suitably defined. Several authors, including Benedetto and Heinig \cite{1}; Heinig and Sinnamon \cite{3}; 
and Jurkat and Sampson \cite{5} have made sizable inroads, but no general conditions on $u$ and $v$, both necessary and sufficient, are known.

However, if in (\twoweight), we replace the weighted $L^p$ spaces with weighted Lorentz spaces, quite complete results exist, thanks largely to the works of Sinnamon \cite{8} and Benedetto and Heinig \cite{1, p. 18}. The weighted Lorentz spaces $\Lambda_p(w)$ and $\Gamma_p(w)$ are respectively defined to be the set of all nonnegative, measurable functions defined on $[0,\infty)$ such that $\|f\|^p_{\Lambda_p(w)}:=\int_0^\infty f^{*p}w < \infty$ and $\|f\|^p_{\Gamma_p(w)}:=\int_0^\infty f^{**p}w < \infty$ where $f^{**}(x)=\frac{1}{x}\int_0^x f^*$. Taking $f \in L^1\cap L^2$ with $p,q\in (0,\infty)$, Sinnamon \cite{8} found necessary conditions and sufficient conditions on nonnegative $u(t)$ and $v(t)$ such that 

$$\|\widehat f\|_{\Lambda_q(u)} \leq C\|f\|_{\Gamma_p(t^{p-2}v(1/t))}\tag\eqnam\Sinnamons\neweq$$

\noindent by exploiting the unweighted version of this inequality due to Jodeit and Torchinsky \cite{4, Theorem 4.6}. When $q=2$ and $0<p\leq 2$ the conditions that Sinnamon gives are {\it both} necessary and sufficient \cite{8, section 5}.  Benedetto and Heinig \cite{1, p. 18} found necessary and sufficient conditions on $u$ and $v$ such that $$\|\widehat f\|_{\Lambda_q(u)} \leq C\|f\|_{\Lambda_p(v)}.\tag\eqnam\Benny\neweq$$

Our corollary below shows that for decreasing functions, the boundedness of a Hardy-type operator implies the two weight Lebesgue inequality (\twoweight). 

\proclaim{Corollary} Let $p,q \in (0,\infty)$. Let $f\in L^p(v)$ and suppose $f\in L^1[0,\infty)$ is nonnegative and decreasing. If there exists a constant $C$ such that the weighted inequality for the Hardy operator

$$\left(\int_0^\infty \left(\int_0^z f(x)\,dx\right)^q \frac{u(1/z)}{z^2}\,dz)\right)^{1/q}\leq C \left(\int_0^\infty f(x)^p v(x)\,dx\right)^{1/p}\tag \eqnam\Hardy\neweq$$

\noindent holds, then

$$\left(\int_0^\infty |\widehat f(z)|^q u(z)\,dz\right)^{1/q} \leq C\left(\int_0^\infty f(x)^pv(x)\,dx.\right)^{1/p}.$$

\noindent That is, there exists a constant $C$ such $\|\widehat f\|_{L^q(u)} \leq C\|f\|_{\Lambda_p(v)}$ for $f$ decreasing.

\endproclaim

Necessary and sufficient conditions on $u$ and $v$ such that (\Hardy) holds are well known, both in the case of general functions $f$ as well as for decreasing functions $f$. The case of general functions is due to the work of many authors, one can consult Maz'ja \cite{6} or Benedetto and Heinig \cite{1, p. 6} as references. Sawyer discovered necessary and sufficient conditions such that (7) holds for decreasing functions \cite{7, Theorem 2}.

Although the corollary only applies to decreasing functions $f$, it has the advantage of estimating the  $L^q(u)$ norm of the Fourier transform as opposed to the norm of the decreasing rearrangement of the Fourier transform as in (\Sinnamons) and (\Benny). These are, in general, not comparable. This is simply because the decreasing rearrangement defined with respect to Lebesgue measure and the weight function $u$ are incompatible. For example, if we take the function $u(x)=\chi_{(1,\infty)}(x)$ and 

$$g(x) = \cases x, & \hbox{for\ } x\in [0,1];\cr
2-x, & \hbox{for\ }x\in (1,2];\cr
0,& \hbox{otherwise,\ }
\endcases$$

\noindent so that $g^*(x)=-.5x+1$, then $\|g\|_{\Lambda_p(u)} = \|g^*\|_{L^p(u)} < \|g\|_{L^p(u)}$. 
The corollary provides us with an estimate for $\|\widehat f\|_{L^q(u)}$ which may in fact be larger than $\|\widehat f\|_{\Lambda_p(u)}$.

\demo{Proof of Corollary} By (\basicfourier) and the fact that $f$ is decreasing we have for $z>0$

$$|\widehat f(z)| \leq C\int_0^{1/z} f(x)\,dx.$$

\noindent Hence, by changing variables and applying the assumption we have

$$\aligned
\left(\int_0^\infty |\widehat f(z)|^qu(z)\,dz\right)^{1/q} &\leq C\left(\int_0^\infty \left(\int_0^{1/z} f(x)\,dx\right)^qu(z)\,dz\right)^{1/q}\cr
&=C\left(\int_0^\infty \left(\int_0^z f(x)\,dx\right)^q\frac{u(1/z)}{z^2}\,dz\right)^{1/q}.\cr
&\leq C\left(\int_0^\infty f(x)^pv(x)\,dx\right)^{1/p}.\cr
\endaligned$$
\enddemo
\endexample

\remark{Acknowledgements} The author would like to thank Dr. Greg Oman for fruitful discussions regarding some topics in this paper. Also the author would like to thank the referee of an earlier version of this work|his/her comments have greatly improved the clarity of this paper.
\endremark

\refstyle{C}

\enddocument